\documentclass[12pt,reqno]{amsart}
\usepackage{amssymb,delarray}
\usepackage{amsfonts}
\usepackage{epsfig}
\usepackage[all]{xy}

\def\leq{\leqslant}
\def\geq{\geqslant}

\newtheorem{thm}{Theorem}
\newtheorem{lem}
{Lemma}
{Proposition}
\newtheorem{claim}
{Claim}

\newtheorem{cor}
{Corollary}
{Remark}

{\catcode`\@=11
\gdef\n@te#1#2{\leavevmode\vadjust{%
 {\setbox\z@\hbox to\z@{\strut#1}%
  \setbox\z@\hbox{\raise\dp\strutbox\box\z@}\ht\z@=\z@\dp\z@=\z@%
  #2\box\z@}}}
\gdef\leftnote#1{\n@te{\hss#1\quad}{}}
\gdef\rightnote#1{\n@te{\quad\kern-\leftskip#1\hss}{\moveright\hsize}}
\gdef\?{\FN@\qumark}
\gdef\qumark{\ifx\next"\DN@"##1"{\leftnote{\rm##1}}\else
 \DN@{\leftnote{\rm??}}\fi{\rm??}\next@}}

\begin{document}
\baselineskip=13.7pt plus 2pt 

\title[Factorization semigroups. II] {Factorization semigroups and irreducible components of
Hurwitz space. II}
\author[Vik.S. Kulikov]{Vik.S. Kulikov}

\address{Steklov Mathematical Institute}
 \email{kulikov@mi.ras.ru}
\dedicatory{} \subjclass{}
\thanks{This
research was partially supported by grants of NSh-4713.2010.1, RFBR
11-01-00185, and by AG Laboratory HSE, RF government grant, ag.
11.G34.31.0023. }
\begin{abstract}
This article is a continuation of the article \cite{K1}. Let {\rm
$\text{HUR}_{d,t}^{\mathcal S_d}(\mathbb P^1)$} be the Hurwitz space
of degree $d$ coverings of the projective line $\mathbb P^1$ with
Galois group $\mathcal S_d$ and having fixed monodromy type $t$
consisting of a collection of local monodromy types (that is, a
collection of conjugacy classes of permutations $\sigma$ of the
symmetric group $\mathcal S_d$ acting on the set
$I_d=\{1,\dots,d\}$). We prove that if the type $t$ contains big
enough number of local monodromies belonging to the conjugacy class
$C$ of an odd permutation $\sigma$ which leaves fixed $f_C\geq 2$
elements of $I_d$, then the Hurwitz space {\rm
$\text{HUR}_{d,t}^{\mathcal S_d}(\mathbb P^1)$} is irreducible.
\end{abstract}

\maketitle
\setcounter{tocdepth}{1}


\def\st{{\sf st}}



\section*{Introduction} \label{introduc}
This article is a continuation of article \cite{K1}. To formulate
the results of presented article, let us recall main definitions and
notations used in \cite{K1}. A collection $(S,G,\alpha,\rho) $,
where $S$ is a semigroup, $G$ is a group, and $\alpha :S\to G$,
$\rho :G\to \text{Aut}(S)$ are homomorphisms, is called {\it a
semigroup $S$ over a group $G$} if for all $s_1,s_2\in S$ we have
\begin{equation} \label{relx}
s_1\cdot s_2=\rho (\alpha (s_1))(s_2)\cdot s_1= s_2\cdot\lambda
(\alpha (s_2))(s_1),\end{equation} where $\lambda(g)=\rho (g^{-1})$.
Let $(S_1,G,\alpha_1,\rho_1)$ and $(S_2,G,\alpha_2,\rho_2)$ be two
semigroups over a group $G$. A homomorphism of semigroups $\varphi :
S_1\to S_2$ is said {\it to be defined over} $G$ if
$\alpha_1(s)=\alpha_2(\varphi(s))$ and
$\rho_2(g)(\varphi(s))=\varphi(\rho_1(g)(s))$ for all $s\in S_1$ and
$g\in G$.

A pair $(G,O)$, where  $O\subset G$ is a subset of a group $G$
invariant under the inner automorphisms, is called an {\it equipped
group}. To each equipped group $(G,O)$ one can associate a semigroup
$S_O=S(G,O)$ over $G$ (called the {\it factorization semigroup of
the elements of $G$ with factors in $O$}) generated by the elements
of the alphabet $X=X_O=\{ x_g \mid g\in O\}$ being subject to the
following relations:
\begin{equation} \label{rel1} x_{g_1}\cdot
x_{g_2}=x_{g_2}\cdot\,x_{g_2^{-1}g_1g_2}= x_{g_1g_2g_1^{-1}}\cdot
x_{g_1}
\end{equation}
for each $x_{g_1}, x_{g_2}\in X$ and if $g_2=\bold 1$ then
$x_{g_1}\cdot x_{\bold 1}=x_{g_1}$. The map $\alpha : X\to G$, given
by $\alpha (x_g)=g$ for each $x_g\in X$, induces a homomorphism
$\alpha :S_O\to G$ called the {\it product homomorphism}. The action
(from the left) $\rho$ of the group $G$ on $S_O$ is defined by the
action on the alphabet $X$ as follows: $$ x_a\in X\mapsto \rho
(g)(x_a)=x_{gag^{-1}}\in X$$ for each $g\in G$. Note that $\alpha
(\rho (g)(s))=g\alpha(s)g^{-1}$ for all $s\in S_O$ and all $g\in G$.

Let $O\setminus \{ {\bf{1}} \}=C_1\sqcup \dots \sqcup C_m$ be the
decomposition of the set $O$ into the disjoint union of conjugacy
classes of elements of $G$. Each element $s=x_{g_1}\cdot \, .\, .\,
.\, \cdot x_{g_n}\in S_O$ defines an element $\tau(s)=n_{1}C_1+\dots
+n_{m}C_m$ of the free abelian semigroup generated by $C_1,\dots,
C_m$ ($\tau(s)$ is called the {\it type} of $s$), where $n_{i}$ is
equal to the number of factors $x_{g_j}$ entering in the
factorization $s=x_{g_1}\cdot \, .\, .\, .\, \cdot x_{g_n}$ for
which $g_j\in C_i$. The total number $n=\sum_{i=1}^m n_{i}$ is
called the length of $s$ and it is denoted by $ln(s)$. A
subsemigroup $S$ of $S_G$ is called {\it stable} if there is an
element $s\in S$ (called a {\it stabilizing element} of $S$) such
that $s_1\cdot s=s_2\cdot s$ for any $s_1,s_2\in S$ such that
$\alpha(s_1)=\alpha(s_2)$ and $\tau(s_1)=\tau(s_2)$.

For an element $s=x_{g_1}\cdot \, .\, .\, .\, \cdot x_{g_n}\in S_O$
one can associate a subgroup $G_s=\langle g_1, \dots ,g_n\rangle$ of
$G$ generated by the elements $g_1,\dots, g_n$. For each two (not
necessary proper) subgroups $H$ and $\Gamma$ of $G$ one can define
subsemigroups $S_{O}^H=\{ s\in S(G,O)\mid G_s=H\}$ and
$S_{O,\Gamma}=\{ s\in S(G,O)\mid \alpha(s)\in \Gamma\}$. If $H$ and
$\Gamma$ are normal subgroups of $G$, then $S_{O,\Gamma}$ and
$S_O^H$ are semigroups over $G$. By definition,
$S_{O,\Gamma}^H=S_{O,\Gamma}\cap S_O^H$.

Let $\mathcal S_d$ be the symmetric group acting on the set $I_d=\{
1,\dots, d\}$ and $T_d\subset \mathcal S_d$ be the subset of
transpositions. The semigroup $S_{\mathcal S_d}$ is denoted by
$\Sigma_d$. By Theorem 2.3. in \cite{K1}, the element
$$h=(\prod_{i=1}^{d-1}x_{(i,i+1)})^3$$
is a stabilizing element of $\Sigma_d$, where $(i,i+1)\in T_d$ is a
transposition permuting the elements $i$ and $i+1$ of $I_d$.

The aim of this article is to generalize this result to the case of
almost all odd elements of the symmetric group $\mathcal S_d$. More
precisely, let $C=C_{\sigma}$ be the conjugacy class of a
permutation $\sigma\in \mathcal S_d$.  Denote by $n_C$ the order of
$\sigma\in C$, by $k_C=|C|$ the number of elements of $C$, and by
$f_C$ the number of elements of $I_d$ fixed under the action of
$\sigma\in C$ on $I_d$.

As is known, if $\sigma$ is an odd permutation, then the elements of
$C$ generate the group $\mathcal S_d$ and, in particular, any
transposition $(i,j)\in \mathcal S_d$ is a product of some
permutations belonging to $C$. In the case  $f_C\geq 2$, denote by
$m_C$ the minimal number (counted with multiplicities) of
permutations of $C\cap \mathcal S_{d-2}$ needed to express $(1,2)$
as a product of permutations of $C$ and fix one of the such
expressions:
\begin{equation} \label{exp} (1,2)=\sigma_1\dots \sigma_{m_C},\qquad \sigma_i\in
C\cap \mathcal S_{d-2}.\end{equation}

\begin{thm} \label{f1} Let $C$ be the conjugacy class of an odd
permutation $\sigma\in \mathcal S_d$. If $f_C\geq 2$ then there is a
constant $$N=N_C<3^{d-3}(2d-1)(d-1)m_C+n_Ck_C+1$$ such that any
element $s=\widetilde s\cdot \overline s\in \Sigma_d^{\mathcal
S_d}$, where  $\overline s\in S_{C}$, is uniquely defined by
$\tau(s)$ and $\alpha(s)$ if $ln(\overline s)\geq N$.
\end{thm}

\begin{cor} Let an equipped symmetric group $(\mathcal S_d, O)$ be
such that the set $O$ contains a conjugacy class $C$ of an odd
permutation $\sigma$, $f_C\geq 2$. Then $S_O=S(\mathcal S_d,O)$ is a
stable semigroup.
\end{cor}

Note that in general case the constant $N_C$, existence of which is
claimed in Theorem \ref{f1}, is greater than $1$. For example, as it
was shown in \cite{W}, this is the case when $C$ is the conjugacy
class of the permutation $\sigma =(1,2)(3,4,5)\in \mathcal S_8$.

The proof of Theorem \ref{f1} is similar to the proof of Theorem 2.3
in \cite{K1} and it is based on the following theorem.

\begin{thm} \label{iso} Let $C$ be the conjugacy class of an odd
permutation $\sigma\in \mathcal S_d$ and an element $\overline
s_{(i_1,i_2)} \in S_C$ be such that
\begin{itemize}
\item[$(i)$] $\alpha(\overline
s_{(i_1,i_2)})=(i_1,i_2)$,
\item[$(ii)$] there are $i_3,i_4\in I_d\setminus \{ i_1,i_2\}$ such that
$\rho((i_3,i_4))(\overline s_{(i_1,i_2)})=\overline s_{(i_1,i_2)}$.
\end{itemize}
Then there is an embedding over $\mathcal S_d$ of the semigroup
$S_{T_d}^{\mathcal S_d}$ in $S_C$.
\end{thm}

Let {\rm $\text{HUR}_{d,b}(\mathbb P^1)$} (resp., {\rm
$\text{HUR}_{d,b}^{G}(\mathbb P^1)$}) be the Hurwitz space of
ramified degree $d$ coverings of the projective line $\mathbb P^1$
(defined over $\mathbb C$) branched over $b$ points (and resp., with
Galois group $G$). In \cite{K1}, it was shown that the irreducible
components of {\rm $\text{HUR}_{d,b}(\mathbb P^1)$} are in one to
one correspondence with the orbits of the action of $\mathcal S_d$
by simultaneous conjugations on $\Sigma_{d,\bf{1}, b}=\{ s\in
\Sigma_{d,\bf{1}}\mid ln(s)=b\}$ (that is, of the action defined by
the homomorphism $\rho$) and if $G=\mathcal S_d$, then the
irreducible components of {\rm $\text{HUR}_{d,b}^{\mathcal
S_d}(\mathbb P^1)$} are in one to one correspondence with the
elements of $\Sigma_{d,{\bf{1}}}^{\mathcal S_d}$ of length equal to
$b$. If an irreducible component of {\rm $\text{HUR}_{d,b}^{\mathcal
S_d}(\mathbb P^1)$} corresponds to an element $s\in
\Sigma_{d,{\bf{1}}}^{\mathcal S_d}$, then we call $\tau(s)$ {\it
monodromy factorization type} of the  coverings belonging to this
component. Denote by {\rm $\text{HUR}_{d,t}^{\mathcal S_d}(\mathbb
P^1)$} the union of irreducible components corresponding to the
elements $s\in \Sigma_{d,{\bf{1}}}^{\mathcal S_d}$ with $\tau(s)=t$.

As a corollary of Theorem \ref{f1} we obtain
\begin{thm} \label{TH1} Let $C$ be the conjugacy class of an odd
permutation $\sigma\in \mathcal S_d$ such that $f_C\geq 2$. If the
monodromy factorization type $t$ contains more than $N_C$ factors
belonging to $C$, where $N_C$ is defined in Theorem {\rm \ref{f1}},
then the space {\rm $\text{HUR}_{d,t}^{\mathcal S_d}(\mathbb P^1)$}
is irreducible.
\end{thm}

Note that fot Hurwitz spaces of $d$-sheeted coverings of the disc
ДХЯЙЮ $\Delta=\{ z\in \mathbb C\mid |z|\leq 1\}$ (respectively, of
$d$-sheeted coverings of the affine line $\mathbb C^1$), a
statement, similar to Theorem \ref{TH1}, is also true.

\section{Proof of Theorem \ref{iso}} Without loss of generality, we
can assume that $(i_1,i_2)=(1,2)$ and $(i_3,i_4)=(3,4)$.

For each $(i,j)\in T_d$ let us choose an element $\sigma_{i,j}\in
\mathcal S_d$ such that $(i,j)=\sigma_{i,j}(1,2)\sigma_{i,j}^{-1}$
and put
$$c=\overline s_{(1,2)}^2\cdot \overline s_{(2,3)}^2\cdot \, .\, .\, .\cdot \overline s_{(d-1,d)}^2,$$ where $\overline
s_{(i,j)}=\rho(\sigma_{i,j})(\overline s_{(1,2)})$.

Obviously, we have $\alpha(\overline s_{(i,j)})=(i,j)$ and
$\alpha(c)=\bf{1}$. Since the transpositions $(1,2),\dots, (d-1,d)$
generate the group $\mathcal S_d$, then $c\in S_{C,\bf{1}}^{\mathcal
S_d}$. Therefore, by Proposition 1.1 (2) in \cite{K1}, the element
$c$ is fixed under the conjugation action of $\mathcal S_d$ on
$S_C$.

For $k\geq 4$ let us denote by $Z_k\simeq \mathcal S_2\times
\mathcal S_{k-2}$ a subgroup of $\mathcal S_d$ generated by
transpositions $(1,2)$ and $(i,j)$, $3\leq i<j\leq k$. Note that
$Z_d$ is the centralizer in $\mathcal S_d$ of the transposition
$(1,2)$.

\begin{claim} \label{cl1} There is an element $z_{(1,2)}\in S_C$ such that
$\alpha(z_{(1,2)})=(1,2)$ and $\rho(\sigma)(z_{(1,2)})=z_{(1,2)}$
for each $\sigma \in Z_d$
\end{claim}
\proof  By induction in $k$, let us show that there is an element
$y_{(1,2),k}\in S_C^{\mathcal S_d}$ such that
$\alpha(y_{(1,2),k})=(1,2)$ and
$\rho(\sigma)(y_{(1,2),k})=y_{(1,2),k}$ for each $\sigma\in Z_k$.
Then $z_{(1,2)}=y_{(1,2),d}$ is a desired element.

Put $y_{(1,2),4}=\overline s_{(1,2)}\cdot c$. If we move the first
factor $\overline s_{(1,2)}$ to the right, then we obtain
$$\begin{array}{ll}
y_{(1,2),4}= & \overline s_{(1,2)}\cdot \overline s_{(1,2)}\cdot
\overline s_{(1,2)}\cdot \overline s_{(2,3)}^2\cdot \, .\, .\, .
\cdot \overline s_{(d-1,d)}^2=
\\ & \rho((1,2))(\overline s_{(1,2)})\cdot \overline s_{(1,2)}\cdot
\overline s_{(1,2)}\cdot \overline s_{(2,3)}^2\cdot \, .\, .\, .
\cdot \overline s_{(d-1,d)}^2= \\ & \rho((1,2))(\overline
s_{(1,2)})\cdot c=\rho((1,2))(\overline s_{(1,2)})\cdot
\rho((1,2))(c)= \\ & \rho((1,2))(\overline s_{(1,2)}\cdot c)=
\rho((1,2))(y_{(1,2),4}),
\end{array}$$
since $c$ is fixed under the conjugation action of $\mathcal S_d$.

Similarly, by assumption of Theorem \ref{iso},
$$\rho((3,4))(y_{(1,2),4})=\rho((3,4))(\overline s_{(1,2)}\cdot c)=
\rho((3,4))(\overline s_{(1,2)})\cdot \rho((3,4))(c) =\overline
s_{(1,2)}\cdot c=y_{(1,2),4}$$ and hence
$\rho(\sigma)(y_{(1,2),4})=y_{(1,2),4}$ for all $\sigma\in Z_4$.

Assume that for some $k\geq 4$, $k<d$, we constructed an element
$y_{(1,2),k}\in S_C^{\mathcal S_d}$ such that
$\alpha(y_{(1,2),k})=(1,2)$ and
$\rho(\sigma)(y_{(1,2),k})=y_{(1,2),k}$ for all $\sigma\in Z_k$.
Consider an element
$y_{(1,2),k}^{\prime}=\rho((k,k+1))(y_{(1,2),k})$. Obviously,
$y_{(1,2),k}^{\prime}$ belongs to $S_C^{\mathcal S_d}$ and it is
easy to see that $\alpha(y_{(1,2),k}^{\prime})=(1,2)$. Hence, the
element $y_{(1,2),k}\cdot y_{(1,2),k}^{\prime}$ belongs to
$S_{C,\bf{1}}^{\mathcal S_d}$ and therefore it is fixed under the
conjugation action of $\mathcal S_d$. Besides, the element
$y_{(1,2),k}^{\prime}$ is fixed under the action of the group
$Z_k^{\prime}$ generated by transpositions $(i,j)\in Z_{k+1}$,
$i,j\neq k$. Indeed, if $(i,j)\in Z_{k}^{\prime}$ and $i,j\neq k+1$,
then
$$\begin{array}{l} \rho ((i,j))(y_{(1,2),k}^{\prime})=\rho
((i,j))(\rho ((k,k+1))(y_{(1,2),k}))= \\ \rho
((i,j)(k,k+1))(y_{(1,2),k})= \rho ((k,k+1)(i,j))(y_{(1,2),k})=
\\ \rho ((k,k+1))(\rho((i,j))(y_{(1,2),k}))=\rho
((k,k+1))(y_{(1,2),k})=y_{(1,2),k}^{\prime}.\end{array}$$

If $(i,k+1)\in Z_{k}^{\prime}$,  then
$$\begin{array}{l} \rho ((i,k+1))(y_{(1,2),k}^{\prime})=\rho
((i,k+1))(\rho ((k,k+1))(y_{(1,2),k}))= \\ \rho
((i,k+1)(k,k+1))(y_{(1,2),k})= \rho ((k,k+1)(i,k))(y_{(1,2),k})=
\\ \rho ((k,k+1))(\rho((i,k))(y_{(1,2),k}))=\rho
((k,k+1))(y_{(1,2),k})=y_{(1,2),k}^{\prime},\end{array}$$ since
$(i,k)\in Z_k$.

Moreover, the elements $y_{(1,2),k}$, $y_{(1,2),k}^{\prime}$
commute. Indeed,
$$\begin{array}{l} y_{(1,2),k}^{\prime}\cdot
y_{(1,2),k}=
\rho(\alpha(y_{(1,2),k}^{\prime}))(y_{(1,2),k})\cdot y_{(1,2),k}^{\prime}= \\
\rho((1,2))(y_{(1,2),k})\cdot y_{(1,2),k}^{\prime}=y_{(1,2),k}\cdot
y_{(1,2),k}^{\prime}.\end{array}$$

Put $y_{(1,2),k+1}:=y_{(1,2),k}^2\cdot y_{(1,2),k}^{\prime}$. It is
easy to see that $y_{(1,2),k+1}\in S_C^{\mathcal S_d}$ and
$\alpha(y_{(1,2),k+1})=(1,2)$. Let us show that
$\rho(\sigma)(y_{(1,2),k+1})=y_{(1,2),k+1}$ for each $\sigma\in
Z_{k+1}$. First of all note that the group $Z_{k+1}$ is generated by
the elements of the groups $Z_k$ and $Z_k^{\prime}$.

For each $\sigma\in Z_k$ we have
$$\begin{array}{l}
\rho(\sigma)(y_{(1,2),k+1})=\rho(\sigma)(y_{(1,2),k}\cdot
y_{(1,2),k}\cdot y_{(1,2),k}^{\prime})=
\\ \rho(\sigma)(y_{(1,2),k})\cdot \rho(\sigma)(y_{(1,2),k}\cdot
y_{(1,2),k}^{\prime})=y_{(1,2),k}\cdot y_{(1,2),k}\cdot
y_{(1,2),k}^{\prime}, \end{array}$$ since the element
$y_{(1,2),k}\cdot y_{(1,2),k}^{\prime}\in S_{C,\bf{1}}^{\mathcal
S_d}$ is fixed under the conjugation action of $\mathcal S_d$.

Similarly,  for each $\sigma\in Z_k^{\prime}$ we have
$$\begin{array}{l}
\rho(\sigma)(y_{(1,2),k+1})=\rho(\sigma)(y_{(1,2),k}^2\cdot
y_{(1,2),k}^{\prime})=\\
\rho(\sigma)(y_{(1,2),k}^2)\cdot
\rho(\sigma)(y_{(1,2),k}^{\prime})=y_{(1,2),k}^{2}\cdot
y_{(1,2),k}^{\prime}=y_{(1,2),k+1}, \end{array}$$ since the element
$y_{(1,2),k}\cdot y_{(1,2),k}\in S_{C,\bf{1}}^{\mathcal S_d}$ is
fixed under the conjugation action of $\mathcal S_d$. Claim
\ref{cl1} is proved. \qed \\

Consider the orbit $X_{T_{C,d}}$ of the element $z_{(1,2)}$  under
the conjugation action of $\mathcal S_d$ on the semigroup $S_C$,
where $z_{(1,2)}$ is the element constructed in the proof of Claim
\ref{cl1} with the help of the element $\overline s_{(1,2)}$.

\begin{claim} \label{cl2} The map $\overline{\alpha} :X_{T_{C,d}}\to
X_{T_d}=\{ x_{(i,j)}\mid (i,j)\in T_d\}$ given by
$\overline{\alpha}(\rho(\sigma)(z_{(1,2)}))=x_{\sigma(1,2)\sigma^{-1}}$
is one-to-one correspondence.
\end{claim}

\proof The map $\overline{\alpha} :X_{T_{C,d}}\to X_{T_d}$ is
surjective, since for each $(i,j)\in T_d$ there is $\sigma\in
\mathcal S_d$ such that $(i,j)=\sigma(1,2)\sigma^{-1}$ and for this
$\sigma$ we have
$$\alpha(\rho(\sigma)(z_{(1,2)}))=\sigma(1,2)\sigma^{-1}=
(i,j),$$
$$\alpha(\overline{\alpha}(\rho(\sigma)(z_{(1,2)})))=\alpha(x_{\sigma(1,2)\sigma^{-1}})=\sigma(1,2)\sigma^{-1}=(i,j) .$$

The order of the group $Z_d$ is equal to $2(d-2)!$. Therefore, by
Claim \ref{cl1}, the number $|X_{T_{C,d}}|$ of the elements of
$X_{T_{C,d}}$ is not more than
$\frac{d!}{2(d-2)!}=\frac{d(d-1)}{2}=|T_d|$ and
hence $\overline{\alpha} :X_{T_{C,d}}\to X_{T_d}$ is one-to-one. \qed \\

Denote by $z_{(i,j)}$ an element $z\in X_{T_{C,d}}$ such that
$\alpha(z)=(i,j)$ and by $S_{T_{C,d}}$ a subsemigroup of $S_C$
generated by the elements $z_{(i,j)}$, $1\leq i,j\leq d$, $i\neq j$.

\begin{claim} \label{cl3} The
subsemigroup $S_{T_{C,d}}$ of $S_C$ is a semigroup over $\mathcal
S_d$ The elements $z_{(i,j)}\in S_{T_{C,d}}$, $1\leq i,j\leq d$,
$i\neq j$, satisfy the following relations
\begin{equation} \label{trrel}
\begin{array}{l}
z_{(i,j)}=z_{(j,i)} 
\, \, \, \text{for all}\,\, \{ i,j\}_{ord} \subset I_d;
\\
z_{(i_1,i_2)}\cdot z_{(i_1,i_3)}=z_{(i_2,i_3)}\cdot
z_{(i_1,i_2)}=z_{(i_1,i_3)}\cdot z_{(i_2,i_3)} 
\, \, \, \text{for all}\,\, \{i_1,i_2,i_3\}_{ord}\subset I_d;\\
z_{(i_1,i_2)}\cdot z_{(i_3,i_4)}=z_{(i_3,i_4)}\cdot z_{(i_1,i_2)} 
\, \, \, \text{for all}\,\, \{ i_1,i_2,i_3,i_4\}_{ord} \subset I_d.
\end{array}\end{equation}
\end{claim}

\proof It directly follows from constructions of the elements
$z_{(i,j)}$ and Claim 1.1 in \cite{K1}. \qed

\begin{claim} \label{cl4} The map $\overline{\alpha}^{-1} :X_{T_{d}}\to X_{T_{C,d}}$
is extended to surjective homomorphism $\overline{\alpha}^{-1}
:S_{T_{d}}\to S_{T_{C,d}}$ of semigroups over group $\mathcal S_d$.
\end{claim}

\proof Note that if to substitute  $x_{(i,j)}$ instead of
$z_{(i,j)}$ in relations (\ref{trrel}), then we obtain the defining
relations in the semigroup $S_{T_{d}}$. Therefore, it follows from
Claim \ref{cl3} that $\overline{\alpha}^{-1}$ can be extended to a
surjective homomorphism semigroups over group $\mathcal S_d$. \qed \\

For $s\in S_{T_{C,d}}$, where $s$ is a product of $n$ generators
$z_{(i,j)}$ of $S_{T_{C,d}}$, we define the $T$-length of $s$ as
$ln_T(s)=n$. We have $ln(s)=ln_T(\overline{\alpha}^{-1}(s))$ for
$s\in S_{T_{d}}$.

Note that it follows from Claim \ref{cl4} that any statement in
\cite{K1},  in which it is claimed that an element of $S_{T_d}$ can
be represented as the product of some generators $x_{i,j}$, is true
for elements of $S_{T_{C,d}}$ if in the statement we change the
elements $x_{(i,j)}$ by $z_{(i,j)}$ and change the lengths by
$T$-lengths.

Define a subsemigroup $S_{T_{C,d}}^{\mathcal S_d,T}$ of
$S_{T_{C,d}}$ as follows: $S_{T_{C,d}}^{\mathcal
S_d,T}=\overline{\alpha}^{-1}(S_{T_d}^{\mathcal S_d}).$

Theorem \ref{iso} follows from
\begin{claim} \label{clx} The restriction of
$\overline{\alpha}^{-1} :S_{T_{d}}\to S_{T_{C,d}}$ to
$S_{T_d}^{\mathcal S_d}$,
$$\overline{\alpha}^{-1}: S_{T_d}^{\mathcal S_d}\to
S_{T_{C,d}}^{\mathcal S_d,T},$$ is an isomorphism of semigroups over
the group $\mathcal S_d$.
\end{claim}
\proof It follows from Theorem 2.1 in \cite{K1} that the
homomorphism $\overline{\alpha}^{-1}: S_{T_d}^{\mathcal S_d}\to
S_{T_{C,d}}^{\mathcal S_d,T}$ is injective. \qed \\

Note also that  Theorem 2.1 in \cite{K1} and Claim \ref{clx} imply
directly the following
\begin{cor} \label{corx} The elements $s$ of the semigroup $S_{T_{C,d}}^{\mathcal S_d,T}$
are defined uniquely by $\alpha(s)$ and $ln_T(s)$.
\end{cor}

\section{Proof of Theorem \ref{f1}}
Let us consider an element $\overline s_{(1,2)}=x_{\sigma_1}\cdot\,
.\, .\, .\cdot x_{\sigma_{m_C}}$, where $\sigma_1, \dots
,\sigma_{m_C}\in C$ are the factors in factorization (\ref{exp}).

If $f_C\geq 2$ then we can and will assume that all $\sigma_i$
entering into factorization (\ref{exp}) belong to the subgroup
$\mathcal S_d^{\{ 3,4\}}\simeq \mathcal S_{d-2}$ of $\mathcal S_d$
the elements of which leave fixed the elements $3,4\in I_d$.
Therefore the element $\overline s_{(1,2)}= x_{\sigma_1}\cdot .\,
.\, .\cdot x_{\sigma_{m_C}}$ satisfies all conditions of Theorem
\ref{iso} and, consequently, the elements $z_{(i,j)}$, constructed
in section 1 with the help of $\overline s_{(1,2)}=
x_{\sigma_1}\cdot .\, .\, .\cdot x_{\sigma_{m_C}}$, define uniquely
a semisubgroup $S_{T_{C,d}}^{\mathcal S_d,T}$  of $S_C$ isomorphic
to $S_{T_d}$ over the group $\mathcal S_d$.

Note that the length of the element $z_{(1,2)}$, constructed in the
proof of Claim \ref{cl1}, is equal to $ln(z_{(1,2)})=
3^{d-4}(2d-1)m_C$ if we start from the element $\overline
s_{(1,2)}=x_{\sigma_1}\cdot\, .\, .\, .\cdot x_{\sigma_{m_C}}$,
where $\sigma_1, \dots ,\sigma_{m_C}$ are the factors of expression
(\ref{exp}).

Denote by
$$h_C=(z_{(1,2)}\cdot z_{(2,3)}\cdot \, .\, .\, .\cdot
z_{(d-1,d)})^3.$$ Rewrite $h_C$ as a product
$$h_C=x_{\sigma_1}\cdot  .\, .\, .\cdot x_{\sigma_L}, \qquad
\sigma_i\in C\, \, \text{for}\, \, i=1, \dots ,L.$$ It is easy to
see that
$$\ln (h_C)=3^{d-3}(2d-1)(d-1)m_C:=L.$$

To prove Theorem \ref{f1}, we need the following

\begin{claim} \label{cl5} Under the conditions of Theorem \ref{f1}, let an element $s =\widetilde s\cdot
\overline s\in \Sigma_d^{\mathcal S_d}$ be such that $\overline s\in
S_C$ of length $ln(\overline s):=M\geq
3^{d-3}(2d-1)(d-1)m_C+n_Ck_C$. Then the element $s$ can be
represented as $s=\widetilde s'\cdot h_C$.
\end{claim}

\proof Let \begin{equation} \label{fuc} \overline
s=x_{\sigma_1}\cdot .\, .\, .\cdot x_{\sigma_M},\end{equation}
$\sigma_i\in C$. Since  $M=ln(\overline s)\geq
3^{d-3}(2d-1)(d-1)m_C+n_Ck_C>n_Ck_C$, then for some $\sigma\in C$
there are at least $n_C+1$ factors entering into factorization
(\ref{fuc}) equal to $\sigma$. Therefore $\overline s$ can be
written in the form: $\overline s=\overline s'\cdot
x_{\sigma}^{n_C}$, where $\overline s'\in S_{C}$ is such that
$\widetilde s\cdot\overline s'\in \Sigma_{d}^{\mathcal S_d}$. By
Lemma 1.1 in \cite{K1}, we have $$s=\widetilde s\cdot \overline
s'\cdot x_{\sigma}^{n_C}=\widetilde s\cdot \overline s'\cdot
x_{\sigma_L}^{n_C}=\widetilde s\cdot \overline s_L\cdot
x_{\sigma_L},$$ where $\overline s_L=\overline s'\cdot
x_{\sigma_L}^{n_C-1}$. Note that $\widetilde s\cdot \overline s_L\in
\Sigma_d^{\mathcal S_d}$ and $ln(\overline s_L)>n_Ck_C$. Therefore,
by the same arguments, the element $\widetilde s\cdot \overline s_L$
can be written in the form: $\widetilde s\cdot \overline
s_L=\widetilde s\cdot \overline s_L'\cdot
x_{\sigma_{L-1}}^{n_C-1}\cdot x_{\sigma_{L-1}}$. Put $\overline
s_{L-1}=\overline s_L'\cdot x_{\sigma_{L-1}}^{n_C-1}$. Repeating the
same arguments for $\widetilde s\cdot \overline s_{L-1}$ we obtain
that $\widetilde s\cdot \overline s_{L-1}=\widetilde s\cdot
\overline s_{L-2}\cdot x_{\sigma_{L-1}}$, and so on. Finally, on the
$L$th step we obtain that $$s=\widetilde s\cdot \overline
s=\widetilde s\cdot \overline s_0\cdot (x_{\sigma_1}\cdot\, .\, .\,
.\cdot x_{\sigma_L})=\widetilde s\cdot \overline s_0\cdot h_C.
\qed$$

Now to complete the proof of Theorem \ref{f1}, recall that the proof
of Theorem 3.2 in \cite {K1} consists of two parts. In the first
part of the proof, for any element $s=\widetilde s\cdot \overline
s\in \Sigma_d^{\mathcal S_d}$, where $\overline s\in S_{T_d}$ has
the length $ln(\overline s)\geq 3(d-1)$, it was proved the existence
of another factorization $s=\widetilde s_1\cdot \overline s_1$ such
that $\overline s_1\in S_{T_d}^{\mathcal S_d}$ with $ln(\overline
s_1)=3(d-1)$. In this case the element $\overline s_1$ is uniquely
determined by its product $\alpha (\overline s_1)=\alpha (\widetilde
s_1)^{-1}\alpha(s)$.

In the second part of the proof of Theorem 3.2 in \cite{K1} it was
proved that for a such factorization $s=\widetilde s_1\cdot
\overline s_1$ there is another factorization $s=\widetilde s_2\cdot
\overline s_2$, where again  $\overline s_2\in S_{T_d}^{\mathcal
S_d}$ has the length $ln(\overline s_2)=3(d-1)$ and $\widetilde s_2$
is uniquely determined by the type $\tau(\widetilde s_1)$. The proof
of the last statement used only properties of the semigroup
$S_{T_d}$ and relations (\ref{relx}) in the factorization
semigroups.  Therefore, by Claims \ref{clx} and \ref{cl5}, the end
of the proof of Theorem \ref{f1} coincides  ОНВРХ ДНЯКНБМН with the
second part of the proof of Theorem 2.3 in \cite{K1}. Only, we must
do the following changes: the elements $x_{(i,j)}$ are changed by
$z_{(i,j)}$, the lengths of elements are changed by $T$-lengths, the
element $h_{d,g}$ is changed by $\overline \alpha^{-1}(h_{d,g})$,
the semigroup $S_{T_d}^{\mathcal S_d}$ is changed by
$S_{T_{C,d}}^{\mathcal S_d,T}$, and the homomorphism $r$ is changed
by $\overline \alpha^{-1} \circ r$. \qed

However, according the request of the referee, this proof is given
once again.  For this purpose, denote by $h_{C,d,g}=\overline
\alpha^{-1}(h_{d,g})$ the image of the Hurwitz element
$h_{d,g}=x_{(1,2)}^{2g} \cdot x_{(1,,2)}^2\cdot .\, .\, .\cdot
x_{(d-1,d)}^2$.

\begin{lem} \label{hhhh1} For any disjoint union   $\{ i_{1,1},\dots, i_{k_1,1}\}
\sqcup \dots \sqcup \{ i_{1,n},\dots, i_{k_n,n}\}$ of ordered
subsets of $I_d$ the Hurwitz element $h_{C,d,0}$ can be represented
as a product
$$ h_{C,d,0}=(z_{(i_{1,1},i_{2,1})}\cdot \,. \, .\, .\, \cdot z_{(
i_{k_1-1,1},i_{k_1,1})})\cdot \, .\, .\, .\, \cdot
(z_{(i_{1,n},i_{2,n})}\cdot \,. \, .\, .\, \cdot z_{(
i_{k_n-1,n},i_{k_n,n})})\cdot \overline h,$$ where $\overline h$ is
an element of $S_{T_{C,d}}^{\mathcal S_d,T}$.
\end{lem}
\proof It directly follows from Lemma 2.9 in \cite{K1} and Claim \ref{clx}. \qed \\

By Claim \ref{cl5}, the element $s$ can be represented as a product:
$s=\widetilde s'\cdot \overline s$, where $\overline s$ is an
element of $S_{T_{C,d}}^{\mathcal S_d,T}$ of $T$-length  $k\geq
3(d-1)$ (in our case $\overline s=h_{C}$ and $k=3(d-1)$), and let
$\widetilde s'=x_{\sigma'_{1}}\cdot \, .\, .\, .\, \cdot
x_{\sigma'_{m}}$. By Proposition 2.4 in \cite{K1} and Claim
\ref{clx}, we have $\overline s=h_{C,d,0}\cdot \overline s'$.

To complete the proof of Theorem \ref{f1}, let us use induction on
$m$. If $m=0$ (that is, if $s\in S_{T_{C,d}}$), then Theorem
\ref{f1} follows from Proposition 2.4 in \cite{K1} and Claim
\ref{clx}.

Let $m=1$. 
For the canonical representative  $\sigma_{m,0}$ of type
$t(\sigma_{m})$ (the definition of the canonical representative is
given in \cite{K1}), there is an element $\overline{\sigma}_m\in
\mathcal S_d$ such that
$\sigma_{m,0}=\overline{\sigma}_m^{-1}\sigma_{m}'\overline{\sigma}_m$.
The permutation $\overline{\sigma}_m$  can be factorized into the
product of cyclic permutations and each cyclic permutation can be
factorized into the product of transpositions:
$$\overline{\sigma}_m=((i_{1,1},i_{2,1})\dots (
i_{k_1-1,1},i_{k_1,1}))\dots ((i_{1,n},i_{2,n})\dots (
i_{k_n-1,n},i_{k_n,n})).$$ Consider the element
$$\overline r(x_{\overline{\sigma}_m})= (z_{(i_{1,1},i_{2,1})}\cdot \,. \, .\, .\,
\cdot z_{( i_{k_1-1,1},i_{k_1,1})})\cdot \, .\, .\, .\, \cdot
(z_{(i_{1,n},i_{2,n})}\cdot \,. \, .\, .\, \cdot z_{(
i_{k_n-1,n},i_{k_n,n})})\in S_{T_{C,d}}.$$

By Lemma \ref{hhhh1}, we have
$$h_{C,d,0}=\overline r(x_{\overline{\sigma}_m})\cdot \overline h_m,$$ where $\overline
h_m$ is an element of $S_{T_{C,d}}^{\mathcal S_d,T}$. We have
$$\begin{array}{ll}
s= &  x_{\sigma_m'}\cdot h_{d,0}\cdot \overline s'=
x_{\sigma_m'}\cdot \overline r(x_{\overline{\sigma}_m})\cdot
\overline h_m\cdot \overline s'=\\  &
 \overline r(x_{\overline{\sigma}_m}) \cdot
x_{\sigma_{m,0}}\cdot \overline h_m\cdot \overline s'=
x_{\sigma_{m,0}}\cdot \overline r(x_{\overline{\sigma}_m'}) \cdot
\overline h_m\cdot \overline s',
\end{array}$$
where
$x_{\overline{\sigma}_m'}=\lambda(\sigma_{m,0})(x_{\overline{\sigma}_m})$.
We have $\overline s_1'=\overline r(x_{\overline{\sigma}_m'}) \cdot
\overline h_m\cdot \overline s'\in S_{T_{C,d}}^{\mathcal S_d,T}$ and
$\alpha(\overline s_1')=\sigma_{m,0}^{-1}\alpha(s)$.  Theorem 2.4 in
\cite{K1} and Claim \ref{clx} imply that $\overline s_1'=\overline
r(x_{\sigma})\cdot h_{C,d,g}$, where $\sigma=\alpha(\overline s_1')=
\sigma_{m,0}^{-1}\alpha(s)$ and
$g=\frac{k-ln_t(x_{\sigma})}{2}-d+1$.

Now assume that Theorem \ref{f1} is proved for all $m<m_0$ and
consider an element
$$
s=x_{\sigma_1}\cdot \, .\, .\, .\, \cdot x_{\sigma_{m_0}}\cdot
\overline s_1,$$ where the element $\overline s_1\in
S_{T_{C,d}}^{\mathcal S_d,T}$ has $T$-length equal $k\geq 3(d-1)$.
We have
$$\begin{array}{l}
s=x_{\sigma_1}\cdot \, .\, .\, .\, \cdot x_{\sigma_{m_0}}\cdot
\overline s_1=  x_{\sigma_2'}\cdot \, .\, .\, .\, \cdot
x_{\sigma_{m_0}'}\cdot x_{\sigma_1}\cdot \overline s_1=
\\
x_{\sigma_2'}\cdot \, .\, .\, .\, \cdot x_{\sigma_{m_0}'}\cdot
x_{\sigma_{1,0}}\cdot \overline s_1'=  x_{\sigma_{1,0}}\cdot
x_{\sigma_2''}\cdot\, .\, .\, .\, \cdot x_{\sigma_{m_0}''}\cdot
\overline s_1',
\end{array}$$
where $\sigma_j'=\sigma_1\sigma_j\sigma_1^{-1}$ and
$\sigma_j''=\sigma_{1,0}^{-1}\sigma_j'\sigma_{1,0}$ for $j=2,\dots,
m$, and the element $\overline s_1'\in S_{T_{C,d}}^{\mathcal S_d,T}$
has $T$-length $ln_T(\overline s_1')=k$. It follows from inductive
assumption that
$$s=  x_{\sigma_{1,0}}\cdot
(x_{\sigma_2''}\cdot\, .\, .\, .\, \cdot x_{\sigma_{m_0}''}\cdot
\overline s_1')= x_{\sigma_{1,0}}\cdot (x_{\sigma_{2,0}}\cdot\, .\,
.\, .\, \cdot x_{\sigma_{m_0,0}}\cdot \overline s_1''),$$ where the
element  $\overline s_1''\in S_{T_{C,d}}^{\mathcal S_d,T}$ has
$T$-length $ln_T(\overline s_1'')=k$. By Proposition 2.4 in
\cite{K1} and Claim \ref{clx}, we have $\overline s_1''=\overline
r(x_{\sigma})\cdot h_{C,d,g}$, where $\sigma =\alpha(\overline
s_1'')=(\sigma_{1,0}\dots \sigma_{m,0})^{-1}\alpha(s)$ and
$g=\frac{k-ln_t(x_{\sigma})}{2}-d+1$.
 \qed

 \ifx\undefined\bysame
\newcommand{\bysame}{\leavevmode\hbox to3em{\hrulefill}\,}
\fi

\end{document}